# THE ASYMPTOTIC BEHAVIOR OF DENSITIES RELATED TO THE SUPREMUM OF A STABLE PROCESS

### By R. A. Doney and M. S. Savov

*University of Manchester and Université Pierre et Marie Curie*

If $X$ is a stable process of index $\alpha \in (0,2)$ whose Lévy measure has density $cx^{-\alpha-1}$ on $(0, \infty)$, and $S_1 = \sup_{0 < t \leq 1} X_t$, it is known that $P(S_1 > x) \sim A\alpha^{-1}x^{-\alpha}$ as $x \to \infty$ and $P(S_1 \leq x) \sim B\alpha\rho^{-1}x^{\alpha\rho}$ as $x \downarrow 0$. [Here $\rho = P(X_1 > 0)$ and $A$ and $B$ are known constants.] It is also known that $S_1$ has a continuous density, $m$ say. The main point of this note is to show that $m(x) \sim Ax^{-(\alpha+1)}$ as $x \to \infty$ and $m(x) \sim Bx^{\alpha\rho-1}$ as $x \downarrow 0$. Similar results are obtained for related densities.

**1. Introduction and results.** This paper was motivated by the following question which arises in connection with some problems in optimal stopping. Let $X$ be a strictly stable process of index $\alpha \in (0,2)$ which has positive jumps so that its Lévy measure has density

$$(1) \qquad \nu(x) = \begin{cases} c_+ x^{-(\alpha+1)}, & x > 0, \\ c_- |x|^{-(\alpha+1)}, & x < 0, \end{cases}$$

where $c_+ > 0, c_- \geq 0$. Assume also that $X$ is not a subordinator. Let $\tau_x$, for $x > 0$, denote the first passage time of $X$ above level $x$, namely,

$$\tau_x = \inf\{t : X_t > x\}.$$

Then does $\tau_x$ have a density function? If so, how does it behave at zero and infinity?

We prefer to rephrase this in terms of the maximum process defined by $S_t = \sup_{s \leq t} X_s$ which we can do because of the obvious identity

$$P(\tau_x > t) = P(S_t \leq x) = P(S_1 \leq xt^{-\eta}),$$

where $\eta = 1/\alpha$, and we have used the scaling property. It is easy to check that $S_t$ has a continuous density function, $m_t$ say, and by scaling we have

---

$m_t(x) = t^{-\eta} m(x/t^\eta)$ where $m$ stands for $m_1$, so we conclude that $\tau_x$ has a continuous density function $h_x(t)$ given by

$$(2) \qquad\qquad h_x(t) = \eta x t^{-\eta-1} m(xt^{-\eta}).$$

Moreover, we can read off the asymptotic behavior of $h_x$ at infinity and zero from the asymptotic behavior of $m$ at zero and infinity, respectively. It is not difficult to make a conjecture about what this behavior is since it is known from [7], Theorem 4A (see also [5], Proposition 4, page 221) that

$$(3) \qquad\qquad P(S_1 > x) \backsim P(X_1 > x) \backsim \frac{A}{\alpha} x^{-\alpha} \qquad \text{as } x \to \infty,$$

and from [7], Theorem 3A (see also Proposition 2, page 219 in [5]) that

$$(4) \qquad\qquad P(S_1 \le x) \backsim \frac{B}{\alpha\rho} x^{\alpha\rho} \qquad \text{as } x \downarrow 0.$$

Here $\rho = P(X_1 > 0)$, and $A$ and $B$ are explicitly known constants which can be expressed in terms of $\alpha, \rho, c_+$ and $c_-$, and the result (4) is also valid in the spectrally negative case when $c_+ = 0 < c_-, \alpha \in (1, 2)$ and $\alpha\rho = 1$. [Recall that (3) does not hold in this case since $P(X_1 > x)$ is exponentially small at $\infty$.] So the obvious conjecture is that, in these cases,

$$(5) \qquad\qquad m(x) \backsim A x^{-(\alpha+1)} \qquad \text{as } x \to \infty$$

and

$$(6) \qquad\qquad m(x) \backsim B x^{\alpha\rho-1} \qquad \text{as } x \downarrow 0.$$

This question turns out to be closely related to a similar question about the density function $\tilde{p}(x)$, say of $\theta_1^{(1)}$ where $(\theta_s^{(t)}, 0 \le s \le t)$ denotes the stable meander of length $t$; informally, this is the stable process conditioned to stay positive up to time $t$. The following, which is the main result in this paper, confirms that these conjectures are true.

THEOREM 1. *For any strictly stable process $X$ which is such that $|X|$ is not a subordinator, (6) holds and there is a constant $C \in (0, \infty)$ such that*

$$(7) \qquad\qquad \tilde{p}(x) \backsim C x^{\alpha\rho} \qquad \text{as } x \downarrow 0.$$

*If $X$ also has positive jumps, (5) holds and also*

$$(8) \qquad\qquad \tilde{p}(x) \backsim A \rho^{-1} x^{-(\alpha+1)} \qquad \text{as } x \to \infty.$$

REMARK 2. Let $f(x) = f_1(x)$ where $f_t(x)$ denotes the density function of $X_t$; this exists, and is continuous and bounded on $\mathbb{R}$ because the characteristic function of $X_t$ is in $L_1$. A great deal is known about $f$ (see, e.g., [13], pages 87–90). In particular, it is known that $f$ has exactly the same behavior as $x \to \infty$ as that stated for $m$ in (5) which is similar also to that of $\tilde{p}$ in (8). However, its behavior as $x \downarrow 0$ is quite different to those of $m$ and $\tilde{p}$, and, in fact,

$$(9) \qquad \lim_{x \downarrow 0} f(x) = D \in (0, \infty),$$

where the explicit value of $D$ can be read off from (14.30) and (14.33) in [13].

REMARK 3. It should be noted that some of our results are already known in the special case that $\alpha \in (1, 2)$ and $c_- = 0$, that is, the spectrally positive case. Here in [6] the semi-explicit form of the Wiener–Hopf factorization was exploited to show that $m(x) = \sum_1^\infty a_n x^{\alpha n - 2}$ where the $a_n$ are given explicitly; since $\rho = 1 - 1/\alpha$, this confirms (6) in this case. Another result in [6] is an expression for the Fourier transform of $m$, and this is used in [11] to show that in this case (5) also holds.

REMARK 4. If we knew in advance that $m$ was ultimately monotone, as $x \to \infty$ or as $x \to 0$, (5) and (6) would follow immediately from (3) and (4), but as we do not have this information, we have to use a different method. Also, in contrast to the case of $X$, we have no explicit knowledge of the characteristic function of $S$, so our method has to be rather indirect.

REMARK 5. We can read off from these results that the asymptotic behavior of the density of $\tau_x$ is given, for each fixed $x \in (0, \infty)$, by

$$(10) \qquad \begin{aligned} h_x(t) &\backsim \eta B x^{\alpha \rho} t^{-(\rho+1)} &&\text{as } t \to \infty \quad \text{and} \\ h_x(t) &\to A\eta/x^\alpha &&\text{as } t \downarrow 0. \end{aligned}$$

**2. Preliminaries on Lévy processes.** In our proofs we will use several identities for stable processes, all of which are in fact special cases of results valid for general Lévy processes. It seems worthwhile to state these results in the general case. Therefore in this section alone, $X$ will denote a generic Lévy process which is not compound Poisson. $S$ and $I$ will be the associated supremum and infimum processes, we will write $L$ and $L^*$ for the local times at zero of the reflected processes $S - X$ and $X - I$, respectively, and $n$ and $n^*$ will denote the characteristic measures of the excursions away from 0 of these processes. We write $\varepsilon$ for a typical excursion, $\zeta$ for its lifetime and $\overline{\pi^*}(t) = n^*(\zeta > t)$. Then

$$\tilde{P}_t(dx) := n^*(\varepsilon_t \in dx | \zeta > t) = \frac{n^*(\varepsilon_t \in dx, \zeta > t)}{\overline{\pi^*}(t)}$$

is a probability distribution which, in the stable case, coincides with that of $\theta_t^{(t)}$, the stable meander of length $t$, at time $t$. (See [5], page 234.) Our first result connects these quantities to the distribution of $S_t$.

Lemma 6. *There is a constant, $0 < k < \infty$, which depends only on the normalization of $L$ and $L^*$, such that, for $t, x > 0$,*

$$
\begin{aligned}
kP(S_t > x) &= \int_0^t n^*(\varepsilon_{t-s} > x, \zeta > t - s)(b + \overline{\pi}(s))\,ds \\
&\quad + \Delta n^*(\varepsilon_t > x, \zeta > t),
\end{aligned}
\tag{11}
$$

*where $\Delta$ denotes the drift, $b$ the killing rate and $\overline{\pi}(s) = n(\zeta > s)$ the tail of the Lévy measure of the increasing ladder time process $T$.*

Proof. First we note that $P(S_t > x) = P(X_t - I_t > x)$ by duality. Next, if $\varepsilon^{(r)}$ denotes the excursion of $X - I$ starting at time $r \in G := \{r : (X - I)_r = 0\}$, we see that

$$
P(X_t - I_t > x) = E(\mathbf{1}_{\{\varepsilon_{t-g_t}^{(g_t)} > x\}}),
$$

where $g_t := \sup(r \le t : r \in G)$. But the compensation formula gives

$$
E(\mathbf{1}_{\{\varepsilon_{(t-g_t)}^{g_t} > x\}}) = E\left(\int_0^t dL^*(s) n^*(\varepsilon_{t-s} > x, \zeta > t - s)\right)
$$

$$
= \int_0^t n^*(\varepsilon_{t-s} > x, \zeta > t - s) E(dL^*(s))
$$

and we conclude by showing that

$$
kE(dL^*(s)) = b\,ds + \overline{\pi}(s)\,ds + \Delta \delta_0(ds).
\tag{12}
$$

To see this, note that

$$
E\left(\int_0^\infty e^{-qs}\,dL^*(s)\right) = E\left(\int_0^\infty e^{-qL^{*-1}(t)}\,dt\right) = \frac{1}{\Phi^*(q)},
$$

where $\Phi^*$ is the Laplace exponent of the decreasing ladder time process $T^* = L^{*-1}$. But, since the Lévy measure of $T$ is given by $\pi(ds) = n(\zeta \in ds)$,

$$
q \int_0^\infty e^{-qs}(\overline{\pi}(s)\,ds + b\,ds + \Delta \delta_0(ds)) = \int_0^\infty (1 - e^{-qs}) n(\zeta \in ds) + b + \Delta q
$$

$$
= \Phi(q),
$$

the Laplace exponent of $T$. Then (12) follows by Laplace inversion and the identity $\Phi^*(q)\Phi(q) = kq$ (see, e.g., [5], page 166). □

We also need two renewal-type equations which involve $n^*$, but first we recall that if we write

$$
G(dt, dx) = \int_0^\infty P(T_s \in dt, H_s \in dx)\,ds
$$

for the renewal measure in the increasing bivariate ladder process $(T, H)$ of $X$, we have, for $t > 0, x > 0$,

$$(13) \qquad G(dt, dx) = k' n^*(\varepsilon_t \in dx, \zeta > t) \, dt = k' \overline{\pi^*}(t) \tilde{P}_t(dx) \, dt,$$

where $k'$ is a constant which depends only on the choice of the normalization of the local time $L^*$ of $X - I$. This fact comes from Theorem 5 of Alili and Chaumont [2], and the following result is also due to Alili and Chaumont.

LEMMA 7.  *For $t, x > 0$ we have*

$$(14) \qquad tG(dt, dx) = \int_{u=0}^{t} \int_{z=0}^{x} G(du, dz) P_z(X_{t-u} \in dx)$$

*and*

$$(15) \qquad xG(dt, dx) = \int_{u=0}^{t} \int_{z=0}^{x} G(du, dz) \frac{x-z}{t-u} P_z(X_{t-u} \in dx).$$

The first of these statements is the obvious analogue of a result for random walks in [4], and, in fact, (14) is stated and proved in [1], the preliminary version of [2]. However it does not appear in [2] which does, however, contain (15) as (4.9). Unfortunately no proof is given there, but a proof can be found in [3].

## 3. Proof of the small time results.

Throughout this section, $X$ will be a stable process, so that the ladder time processes $T$ and $T^*$ are stable subordinators of index $\rho$ and $(1 - \rho)$, respectively, and we will choose a normalization of the local times so that their exponents are given by $\Phi(q) = q^\rho$ and $\Phi^*(q) = q^{1-\rho}$; this choice entails that $n(\zeta > t) = \overline{\pi}(t) = t^{-\rho}/\Gamma(1 - \rho)$ and $n^*(\zeta > t) = \overline{\pi^*}(t) = t^{\rho-1}/\Gamma(\rho)$. Note also that $\Delta = b = 0$, and the constant $k$ in Lemma 6 is equal to 1. Recall that $\eta = 1/\alpha$, and that the density $f_t$ of $X_t$ is bounded and continuous on $\mathbb{R}$, and has the scaling property

$$(16) \qquad f_t(x) = t^{-\eta} f(xt^{-\eta}) \qquad \text{where } f = f_1.$$

LEMMA 8.  $\theta_t^{(t)}$ *has a continuous and bounded density, $\tilde{p}_t(x)$ say, $G(dt, dx)$ has a continuous bivariate density $g(t, x)$, and these are related by*

$$(17) \qquad g(t, x) = k'' t^{\rho-1} \tilde{p}_t(x) \qquad \text{for } t > 0, x > 0.$$

*Furthermore, $m$ satisfies*

$$(18) \qquad m(x) = \frac{\sin \rho \pi}{\pi} \int_0^1 \frac{\tilde{p}_s(x)}{s^{1-\rho}(1-s)^\rho} \, ds = \frac{\sin \rho \pi}{\pi} \int_0^1 \frac{s^{-\eta} \tilde{p}(xs^{-\eta})}{s^{1-\rho}(1-s)^\rho} \, ds,$$

*where we have written $\tilde{p}$ for $\tilde{p}_1$.*

PROOF.    Note first that if we write $\sigma_0 = \inf\{t : X_t < 0\}$ and introduce the measure, $Q_y(X_t \in dx) := P_y(X_t \in dx, I_t > 0) = P_y(X_t \in dx, \sigma_0 > t)$, then, for $x, y > 0$,

$$P_y(X_t \in dx, \sigma_0 \leq t) = \int_0^t \int_{-\infty}^0 P_y(\sigma_0 \in ds, X_{\sigma_0} \in dz) P_z(X_{t-s} \in dx)$$

and we see that $Q_y(X_t \in \cdot)$ has a continuous and bounded density function given by

$$q_t(y, x) = f_t(x - y) - \tilde{f}_t(y, x),$$

where

$$\tilde{f}_t(y, x) = \int_0^t \int_{-\infty}^0 P_y(\sigma_0 \in ds, X_{\sigma_0} \in dz) f_{t-s}(x - z).$$

Then using the Markov property of the meander and the fact that, given $\varepsilon_t = y$, $\varepsilon_{t+s}$ has law $Q_y(X_s \in \cdot)$, we see that

$$(19) \qquad n^*(\varepsilon_t \in dx, \zeta > t) = \int_0^\infty n^*(\varepsilon_{t/2} \in dy, \zeta > t/2) q_{t/2}(y, x) \, dx$$

and the existence, continuity and boundedness of $\tilde{p}_t(\cdot)$ follows. (This fact is also contained in [12], Theorem 6.) Then (17) follows from (13). Of course $\tilde{p}$ also has the scaling property (16), so specializing Lemma 6 to the stable context and using the identity $\Gamma(\rho)\Gamma(1 - \rho) = \pi/\sin \rho\pi$ gives (18).  □

REMARK 9.    Alternatively, (18) could be deduced from [5], Proposition 16, page 234.

We start by showing that we can deduce the behavior of $\tilde{p}$ at zero from the fact that it is bounded.

PROPOSITION 10.    *(7) holds, viz.* $\tilde{p}(x) \backsim Cx^{\alpha\rho}$ *as* $x \downarrow 0$.

PROOF.    First we write $\kappa(x)$ for $\kappa_1(x)$ with $\kappa_t(x) = \overline{\pi^*}(t)\tilde{p}_t(x)$ the density of $n^*(\varepsilon_t \in dx, \zeta > t)$, and note that it suffices to show that

$$(20) \qquad \lim_{x \downarrow 0} x^{-\alpha\rho} \kappa(x) \in (0, \infty).$$

Next, we have seen that $Q_y(X_t \in dx) = P_y(X_t \in dx, I_t > 0) = q_t(y, x) \, dx$ for $x, y > 0$, and by duality (specifically, [5], Theorem 5, page 47), we deduce that $q_t(y, x) = q_t^*(x, y)$ where

$$q_t^*(x, y) \, dy = P_{-x}(-X_t \in dy, S_t < 0).$$

Using this in the obvious identity

$$\kappa_2(x) = \int_0^\infty \kappa(y) q_1(y, x) \, dy$$

gives

$$\kappa_2(x) = \int_0^\infty \kappa(y) q_1^*(x, y) \, dy \le c \int_0^\infty q_1^*(x, y) \, dy$$

$$= c P_{-x}(S_1 < 0) = c P_0(S_1 < x) \backsim c x^{\alpha\rho}.$$

(Here, and later, $c$ denotes a generic positive constant whose value can change from line to line.) Then, for example, from [8], we recall that the law of "$-X$ starting at $x > 0$ and conditioned to stay positive" has, at time 1, density

$$p_1^{*\uparrow}(x, y) := q_1^*(x, y) \frac{y^{\alpha\rho}}{x^{\alpha\rho}}.$$

Thus if we write $g(y) = \kappa(y) y^{-\alpha\rho}$ which, by the previous result and scaling, is bounded, we see that in the obvious notation,

$$\kappa_2(x) x^{-\alpha\rho} = \int_0^\infty y^{-\alpha\rho} \kappa(y) y^{\alpha\rho} q_1^*(x, y) x^{-\alpha\rho} \, dy$$

$$= \int_0^\infty g(y) p_1^{*\uparrow}(x, y) \, dy$$

$$\to \int_0^\infty g(y) p_1^{*\uparrow}(y) \, dy \in (0, \infty),$$

where $p_1^{*\uparrow}(\cdot)$ denotes the density of 11 "$-X$ starting at 0 and conditioned to stay positive at time 1," and the convergence follows from a result which is stated in [9] and proved in [10]. By scaling, we deduce (20) and then (7). $\square$

REMARK 11. Theorem 6 of [12] gives a completely different proof of Proposition 10, essentially using random walk approximation.

REMARK 12. Since it is known that the density function of $X$, starting from 0 and conditioned to stay positive, is given at time 1 by $p^\uparrow(x) = c x^{\alpha(1-\rho)} \kappa(x)$ (see, e.g., [8]), we deduce from this, and later from Proposition 17, that

(21)    $p_1^\uparrow(x) \backsim c x^\alpha$    as $x \downarrow 0$    and    $p_1^\uparrow(x) \backsim c x^{-(\alpha\rho+1)}$    as $x \to \infty$.

In particular we note that $p_1^\uparrow$ and $p_1^{*\uparrow}$ have the same asymptotic behavior at 0, up to multiplication by a constant.

PROPOSITION 13.    *(6) holds, viz. $m(x) \backsim B x^{\alpha\rho-1}$ as $x \downarrow 0$.*

PROOF.    Write $s^{-\eta} = z$ in (18) and then $zx = y$ to get

$$(22) \quad m(x) = \frac{\alpha \sin \rho\pi}{\pi} \int_1^\infty \frac{\tilde{p}(xz)}{(z^\alpha - 1)^\rho} \, dz = \frac{\alpha x^{\alpha\rho-1} \sin \rho\pi}{\pi} \int_x^\infty \frac{\tilde{p}(y)}{(y^\alpha - x^\alpha)^\rho} \, dy.$$

Since we have seen that $\tilde{p}$ is bounded,

$$\lim_{\delta\downarrow 0} \int_x^{(1+\delta)x} \frac{\tilde{p}(y) \, dy}{(y^\alpha - x^\alpha)^\rho} \leq c \lim_{\delta\downarrow 0} \int_x^{(1+\delta)x} \frac{dy}{(y^\alpha - x^\alpha)^\rho}$$

$$= c x^{1-\alpha\rho} \lim_{\delta\downarrow 0} \int_1^{1+\delta} \frac{dz}{(z^\alpha - 1)^\rho} = 0,$$

as the integral is finite because $\rho < 1$. Moreover on $((1+\delta)x, \infty)$, $(y^\alpha - x^\alpha)^{-\rho} \leq y^{-\alpha\rho}(1 - (1+\delta)^{-\alpha})^{-\rho}$, and, since we know $y^{-\alpha\rho}\tilde{p}(y)$ is integrable on $(0, \infty)$ by Proposition 10, for any $\delta > 0$, dominated convergence gives

$$\int_{(1+\delta)x}^\infty \frac{\tilde{p}(y)}{(y^\alpha - x^\alpha)^\rho} \, dy \to \int_0^\infty y^{-\alpha\rho}\tilde{p}(y) \, dy < \infty.$$

Thus $\lim_{x\downarrow 0} x^{1-\alpha\rho} m(x) \in (0, \infty)$, and since (4) holds, the limit must be $B$. □

REMARK 14.    It is surprising that the result of Proposition 10 only plays a rôle in our proof of Proposition 13 in the case $\alpha\rho = 1$. However it is clear from (18) that the asymptotic behavior of $m$ and $\tilde{p}$ at infinity are closely linked, and we exploit this in the next section.

**4. Proof of the large time results.**    The first step uses the following obvious result:

LEMMA 15.    *Put $\tau = \inf\{t : \varepsilon_t > 1\}$; then for $x > 1$*

$$(23) \quad \kappa(x) = \int_{u=0}^1 \int_{y=1}^\infty n^*(\tau \in du, \varepsilon_\tau \in dy) q_{1-u}(y, x)$$

$$(24) \quad \leq \int_{u=0}^1 \int_{y=1}^\infty n^*(\tau \in du, \varepsilon_\tau \in dy) f_{1-u}(x - y).$$

It does not seem possible to deduce the asymptotic behavior of $\kappa$ from (23), but from (24) we can get a useful upper bound.

PROPOSITION 16.    *We have*

$$(25) \quad \limsup_{x\to\infty} x^{\alpha+1}\kappa(x) < \infty.$$

Proof.  Recall that

$$\lim_{x\to\infty} x^{\alpha+1} f(x) = A. \tag{26}$$

Write the integral in (24) as $I_1 + I_2$, and note that

$$
\begin{aligned}
x^{\alpha+1} I_1 &= \int_{u=0}^{1} \int_{|y-x|>x/2, y>1} (1-u) n^*(\tau \in du, \varepsilon_\tau \in dy) \\
&\qquad\qquad \times \left(\frac{x}{(1-u)^\eta}\right)^{\alpha+1} f\left(\frac{(x-y)}{(1-u)^\eta}\right) \\
&\to A \int_{u=0}^{1} \int_{y>1} (1-u) n^*(\tau \in du, \varepsilon_\tau \in dy) \\
&= A \int_{u=0}^{1} (1-u) n^*(\tau \in du) < \infty,
\end{aligned}
$$

where we have used (26) and dominated convergence. Also when $x$ is sufficiently large,

$$
\begin{aligned}
x^{\alpha+1} I_2 &= \int_{u=0}^{1} \int_{x/2}^{3x/2} (1-u)^{-\eta} x^{\alpha+1} n^*(\tau \in du, \varepsilon_\tau \in dy) f\left(\frac{(x-y)}{(1-u)^\eta}\right) \\
&= \int_{u=0}^{1} \int_{-x/2}^{x/2} (1-u)^{-\eta} x^{\alpha+1} n^*(\tau \in du, \varepsilon_\tau \in x - dz) f\left(\frac{z}{(1-u)^\eta}\right).
\end{aligned}
$$

But if we put $\overline{\nu}(y) = \int_y^\infty \nu(z)\, dz$, then for all sufficiently large $x$,

$$
\begin{aligned}
x^{\alpha+1} \sup_{0<w<1, -x/2<z<x/2} &\ n^*(\varepsilon_\tau \in x - dz \,|\, \varepsilon_{\tau-} = w, \tau = u) \\
&= x^{\alpha+1} \sup_{0<w<1, -x/2<z<x/2} \frac{\nu(x-z-w)\, dz}{\overline{\nu}(1-w)} \le c\, dz.
\end{aligned}
$$

So it follows that

$$
\begin{aligned}
\limsup_{x\to\infty} x^{\alpha+1} I_2 &\le c \int_{u=0}^{1} \int_{-\infty}^{\infty} (1-u)^{-\eta} n^*(\tau \in du) f\left(\frac{z}{(1-u)^\eta}\right) dz \\
&= c \int_{u=0}^{1} \int_{-\infty}^{\infty} n^*(\tau \in du) f(y)\, dy = c n^*(\tau \le 1) < \infty
\end{aligned}
$$

and (25) holds.  □

To exploit this we specialize (15) to get the following integral equation for $g$:

$$x g(t, x) = \int_{u=0}^{t} \int_{z=0}^{x} g(u, z) \frac{x-z}{t-u} f_{t-u}(x-z)\, du\, dz. \tag{27}$$

We are now in a position to establish the behavior of $m(x)$ for large $x$.

PROPOSITION 17. *If $X$ is any strictly stable process which has positive jumps and is not a subordinator, then (5) and (8) hold, viz. $m(x) \backsim \rho \tilde{p}(x) \backsim Ax^{-(\alpha+1)}$ as $x \to \infty$.*

PROOF. Using (17) we can rewrite (27) with $t = 1$ as

$$x^{\alpha+1}\tilde{p}(x) = x^{\alpha+1} \int_{u=0}^{1} \int_{z=0}^{x} u^{\rho-1}\tilde{p}_u(z) \frac{x-z}{x(1-u)} f_{1-u}(x-z) \, du \, dz$$

(28)
$$:= J_1 + J_2,$$

where, by scaling,

$$J_1 = \int_{u=0}^{1} \int_{z=0}^{x/2} u^{\rho-1}\tilde{p}_u(z) \left(1 - \frac{z}{x}\right) \left[\frac{x}{(1-u)^\eta}\right]^{\alpha+1} f((x-z)(1-u)^{-\eta}) \, du \, dz.$$

Using (26) and dominated convergence gives

$$\lim_{x\to\infty} J_1 = A \int_{u=0}^{1} \int_{z=0}^{\infty} u^{\rho-1}\tilde{p}_u(z) \, du \, dz = A/\rho.$$

Using the scaling property of $\tilde{p}$ and a change of variable, we can write

$$J_2 = \frac{1}{x} \int_{u=0}^{1} \int_{y=0}^{x/2} u^\rho \left[\frac{x}{u^\eta}\right]^{\alpha+1} p\left(\frac{x-y}{u^\eta}\right)(1-u)^{-1-\eta} y f(y(1-u)^{-\eta}) \, du \, dy.$$

Now it follows from (25) and the identity, $\kappa(x) = \overline{\pi^*}(1)\tilde{p}(x)$, that the expression $[\frac{x}{u^\eta}]^{\alpha+1}\tilde{p}(\frac{x-y}{u^\eta})$ is bounded above on the range of integration by a constant for all sufficiently large $x$, so that for such $x$,

$$J_2 \leq \frac{c}{x} \int_{u=0}^{1} \int_{y=0}^{x/2} u^\rho (1-u)^{-1-\eta} y f(y(1-u)^{-\eta}) \, du \, dy$$

$$= \frac{c}{x} \int_{u=0}^{1} \int_{z=0}^{x/2(1-u)^\eta} u^\rho (1-u)^{\eta-1} z f(z) \, du \, dz.$$

Since $\int_0^x z f(z) \, dz$ is $O(1), O(\log x)$, or $O(x^{1-\alpha})$ according as $\alpha > 1, \alpha = 1$, or $\alpha < 1$, we check easily that $\lim_{x\to\infty} J_2 = 0$, and hence that $\lim_{x\to\infty} x^{\alpha+1}\tilde{p}(x) = A/\rho$. Putting this into (18) we get

$$\lim_{x\to\infty} x^{\alpha+1}m(x) = \frac{A}{\rho} \frac{\sin\rho\pi}{\pi} \int_0^1 \frac{s^\rho ds}{(1-s)^\rho} = A. \qquad \square$$

REMARK 18. The reader might like to check that the proof above fails if we use the density version of the more obvious identity (14) rather than (15).

**5. Derivatives.** It is not difficult to deduce from the known properties of $f$ that $\tilde{p}$ and $m$ are differentiable $k$ times and that these derivatives satisfy analogues of (18) and (28); it is also known that

$$(29) \qquad f^{(k)}(x) \backsim A_k x^{-(\alpha+k+1)} \qquad \text{as } x \to \infty,$$

where

$$(30) \qquad A_k = \frac{(-)^k \Gamma(\alpha+1+k) A}{\Gamma(\alpha+1)}.$$

This suggests that the derivatives of $m$ and $\tilde{p}$ have asymptotic behaviors which are consistent with the results of Theorem 1. In fact, we have been able to prove such results for large $x$, but have not been able to settle the question for small $x$. (Details of these results can be supplied by the authors on request.)

School of Mathematics
University of Manchester
Manchester M13 9PL
United Kingdom
E-mail: rad@ma.man.ac.uk

Laboratoire de Probabilités
Université Pierre et Marie Curie
Paris 75013
France
E-mail: mladensavov@hotmail.com